\documentclass[11 pt]{amsart}
\usepackage{comment}
\usepackage[pagewise]{lineno}
\usepackage[top=1.4in, bottom=1.25in, left=1.25in, right=1.25in]{geometry}

\usepackage{mathrsfs}

\usepackage{hyperref, amsmath, amsthm, amsfonts, amscd, flafter,epsf,amssymb}

\newtheorem{lemma}{Lemma}
 \newtheorem*{thm 1}{Theorem 1}
 \newtheorem{innercustomthm}{Theorem}
\newenvironment{customthm}[1]
  {\renewcommand\theinnercustomthm{#1}\innercustomthm}
  {\endinnercustomthm}

\newcommand{\be}{\begin{equation}}
\newcommand{\ee}{\end{equation}}

\newcommand{\s}{\sigma}
\newcommand{\di}{i} 

\numberwithin{equation}{section}
\usepackage{parskip}
\begin{document}
\title{Extreme values of Dirichlet $L$-functions at critical points of the zeta function}
\author{Shashank Chorge}
\address{Department of Mathematics, University of Rochester, Rochester NY 14627 USA}
 
\maketitle

\begin{abstract}  
We estimate large  and small values of $|L(\rho',\chi)|$, where $\chi$ is a primitive character mod $q$ for $q>2$ and $\rho'$ runs over  critical points of 
the Riemann  zeta function in the right half of the one-line and on the one-line, that is,  the points where $\zeta'(\rho')=0$ and $1\leq\Re \rho'$.\\
\end{abstract}

    \section{Introduction and previous results}

The Riemann zeta function is defined as $\sum_{n=1}^{\infty}n^{-s}$, where $s$ is a complex number with $\Re s>1$. This function can be meromorphically continued over the whole complex plane. 

Dirichlet $L$-functions are widely studied analogues of the Riemann zeta function. A Dirichlet $L$-function is defined as $\sum_{n=1}^{\infty}\chi(n)n^{-s}$, where $\chi$ is a Dirichlet character mod $q$ for some integer $q$. The Riemann zeta function and Dirichlet $L$-functions play an extremely crucial role in the theory of primes. One of the most famous conjectures, the Riemann hypothesis, states that all the non-trivial zeros of the zeta function lie on the line $\Re s=1/2$. The analogous conjecture for Dirichlet $L$-functions, The generalized Riemann hypothesis, states that all non trivial zeros of such functions lie on  $\Re s=1/2$.

There are many interesting open problems about the Riemann zeta function and Dirichlet $L$-functions. One of those problems is studying extreme values of  Dirichlet $L$-function at the critical points of the zeta function.
 We expect extreme values that an $L$-function takes at critical points of the Riemann zeta function to be very close to the extreme values that the $L$-function would otherwise take to the right of the vertical line $\Re s=1 $. That is, an $L$-function is expected to behave in a manner that is independent of the nature of the points that are special with respect to the Riemann zeta function. The results obtained in this paper corroborate this behavior. In fact, the results in this paper show that extreme values of a Dirichlet $L$-function at critical points of the zeta function are identical to its extreme values at randomly sampled points on $\Re s=1$.

Gonek and Montgomery
~\cite{paper} obtained results sampling zeta at the critical points of the Riemann zeta function to the right of $\Re s=1.$
Let $\rho'=\beta'+\di\gamma'$ denote a typical critical point of the zeta function, that is, a point where $\zeta'(\rho')=0$ and let $C_0$ be Euler's constant.
Assuming RH, they showed that for $\beta' \geq 1$,
\be\label {MG 1}
\limsup\limits_{\gamma' \rightarrow\infty}\frac{|\zeta(\rho')|}{\log\log \gamma'}\leq \frac{1}{2}e^{C_0}
\ee
and, unconditionally, for $\beta'>1$, that
\be\label {MG 2}
\limsup\limits_{\gamma' \rightarrow\infty}\frac{|\zeta(\rho')|}{\log\log \gamma'}\geq \frac{1}{4}e^{C_0}.
\ee

Assuming the Riemann hypothesis (RH), Littlewood \cite{paper3}  proved that
 \be\label {Litt 1}
 \limsup\limits_{t \rightarrow\infty}\frac{|\zeta(1+\di t)|}{\log\log t}\leq 2e^{C_0}.
 \ee

He also proved, unconditionally,  that 
\be\label {Litt 2}
\limsup\limits_{t \rightarrow\infty}\frac{|\zeta(1+\di t)|}{\log\log t}\geq e^{C_0}.
\ee
We can see that the limits obtained in \eqref{MG 1} and \eqref{MG 2} are $1/4$ times the limits obtained in \eqref{Litt 1} and \eqref{Litt 2}.

For small values Gonek and Montgomery proved unconditionally  that 
$$\liminf\limits_{\gamma' \rightarrow\infty}{|\zeta(\rho')|}{\log\log \gamma'}\leq \frac{2\pi^2}{3}e^{-C_0}$$
and assuming RH they proved
$$\liminf\limits_{\gamma' \rightarrow\infty}{|\zeta(\rho')|}{\log\log \gamma'}\geq \frac{\pi^2}{3}e^{-C_0}.$$

In contrast the extreme values of $L$-functions at critical points of the zeta function to the right of $\Re s=1$ are the same as extreme values of the $L$-function on the $\Re s=1$. Results analogous to those of\eqref{Litt 1} and \eqref{Litt 2}  for the small values were proven by Littlewood \cite{paper4} and Titchmarsh \cite{paper5}.

\subsection{More recent results:\\}
The author \cite{papernew5} of this paper extended the result for the extreme values of the zeta function at its critical points in the critical strip. His results are as follows:\\
Let $\sigma_1$ be such that $1/2<\sigma_1<1$. 
Let $\rho'=\beta'+\di \gamma'$ denote a critical point of the Riemann zeta function such that $\sigma_1<\beta'<1$ . Then for any $\epsilon$
and $\epsilon'>0$ and for infinitely many $\rho'$ with $\gamma'\to\infty$, we have unconditionally that
\begin{equation*}\label{eqn thm 2}
 \log|\zeta(\rho')|\geq 
 (B(\sigma_1)  -\epsilon' )\frac{(\log\gamma')^{1-\beta'}}{(\log\log\gamma')^{5-3\beta'+\epsilon}},
  \end{equation*}
  where
  \begin{equation*}
  B(\sigma_1)  =\frac{(\sigma_1-1/2)^{1-\s_1}\log 2}{(1-\s_1)^2 4^{1-\s_1}}.
\end{equation*}

And also:\\
Let $\sigma_1$ be such that $1/2<\sigma_1<1$.
Let $\rho'=\beta'+\di \gamma'$ denote a critical point of the Riemann Zeta function such that $\sigma_1<\beta'<1$. Then for any $\epsilon$ and $\epsilon'>0$ and for infinitely many $\rho'$ with $\gamma'\to\infty$, we have unconditionally that
\begin{equation*}\label{eq1.9}
 \log|\zeta(\rho')|\leq (-B(\s_1)+\epsilon')\frac{(\log\gamma')^{1-\beta'}}{(\log\log\gamma')^{5-3\beta'+\epsilon}}.
\end{equation*}

The recent results by Christoph Aistleitner, Kamalakshya Mahatab and Marc Munsch \cite{papernew1} showed a much improved lower bound on the extreme values of the Riemann zeta function on the line $\Re s=1$. Using resonance method, they showed that there is a constant $C$ such that $$\max{\atop t\in [\sqrt{T},T]}|\zeta(1+\di t)|\geq e^{C_0}(\log\log T+\log\log\log T-C)$$ for all sufficiently large $T$.
Using the same method C. Aistleitner, K. Mahatab, M. Munsch and A. Peyrot \cite{papernew2} proved that for any given $\epsilon$, for all sufficiently large $q$, there is a non principal character $\chi\mod q$ such that
$$|L(1,\chi)|\geq e^{C_0}(\log\log q+\log\log\log q)-C-\epsilon,$$
where $C=1+2\log 2$.
Similar results for the upper bound have been obtained by Y. Lamzouri, Y., X. Li and K. Soundararajan \cite{papernew4}. Assuming GRH they proved that, for any primitive character $\chi \mod q$ with $q\geq 10^{10}$ we have
$$|L(1,\chi)|\leq 2e^{C_0}(\log\log q -\log 2+1/2+1/\log\log q) $$ and
$$\frac{1}{|L(1,\chi)|}\leq \frac{12e^{C_0}}{\pi^2}\bigg(\log\log q -\log 2+1/2+\frac{1}{\log\log q}+\frac{14\log\log q}{\log q}\bigg).$$

\section{Main results of the paper }
The following two theorems give us the bounds for large values of a Dirichlet $L$-function at the critical points of the zeta function.
\begin{customthm}{1}
Assume the Generalized Riemann Hypothesis (GRH). Let $\rho'=\beta'+\di\gamma'$, where $\zeta'(\rho')=0$ and $\beta'\geq 1$. Then for a primitive character $\chi$ mod $q$, 
\begin{equation}\label{th1 est}
  \limsup\limits_{\gamma' \rightarrow\infty}\frac{|L(\rho',\chi)|}{\log\log \gamma'}\leq 2e^{C_0}\frac{\phi(q)}{q},   
    \end{equation}
    where $\phi(q)$ denotes Euler's $\phi$-function.
	\end{customthm}

	\begin{customthm}{2}
	Let $\rho'=\beta'+\di\gamma'$, where $\zeta'(\rho')=0$ and $\beta'> 1$. Then for a primitive character $\chi$ mod $q$,  
\begin{equation}\label{th2 est}
    \limsup\limits_{\gamma' \rightarrow\infty}\frac{|L(\rho',\chi)|}{\log\log \gamma'}\geq e^{C_0}\frac{\phi(q)}{q}.
\end{equation}

	\end{customthm}	

The following two theorems give us the bounds for small values of a Dirichlet $L$-function at the critical points of the zeta function.
	\begin{customthm}{3}
	Assume the GRH. Let $\rho'=\beta'+\di\gamma'$,where $\zeta'(\rho')=0$ and $\beta'\geq 1$. Then for a primitive character $\chi$ mod $q$,   
	\begin{equation}\label{th3 est}
	  \liminf\limits_{\gamma' \rightarrow\infty}|L(\rho',\chi)|\log \log \gamma'\geq\frac{\pi^2e^{-C_0}}{12}\prod_{p\mid q}\frac{p+1}{p}.  
	\end{equation}
	\end{customthm}

	\begin{customthm}{4}
	Let $\rho'=\beta'+\di\gamma'$, where $\zeta'(\rho')=0$ and $\beta'> 1$. Then for a primitive character $\chi$ mod $q$, 
\begin{equation}\label{th4 est}
    \liminf\limits_{\gamma' \rightarrow\infty}|L(\rho',\chi)|\log\log \gamma'\leq \frac{\pi^2e^{-C_0}}{6}\prod_{p\mid q}\frac{p+1}{p}.
\end{equation}	
	
	\end{customthm}

It can be seen using Kronecker's theorem that the bounds in Theorems 1 to 4 are valid, as $t\rightarrow \infty$, even if we replace a critical point $\rho'$ of the Riemann zeta function with a generic point $s=\sigma+\di t$ with $\sigma>1$. This shows that the critical points of the zeta function behave exactly the same way as the randomly sampled points of $L(s,\chi)$ on the line $\Re s=1$.  
The reason we included Theorems 1 and 3 in this paper is to show that the bound in Theorem 1 is twice the bound in Theorem 2 and the bound in Theorem 4 is twice the bound in theorem 3. Thus these bounds essentially follow the similar pattern that can be observed in Littlewood's \cite{paper4} paper and also Montgomery and Gonek's \cite{paper} paper. It is conjectured that the equality should hold in Theorem 2 and 4. Thus it is unlikely that the resonance method will improve the bounds of Theorem 2 and Theorem 4. But it will be interesting to see if the similar bounds can be obtained by the resonance method. 
	
In the remainder of the paper, we use the following notations. The greatest common divisor of $a,b$ will be denoted by $(a,b)$. The letter $p$ denotes a prime and the letters $q$, $n$ represent positive integers.
$\Lambda(n)$ represents the Von-Mangoldt function, which is defined as, 
$$
\Lambda(n):=
\begin{cases}
&\log p, \text{\qquad if $n=p^k$ for some prime $p$ and positive integer $k$}\\
           &  0,     \text{\qquad\qquad    otherwise}.
\end{cases}
$$

\section{\textbf{Lemmas for the proof of Theorem 1}}
In this section we prove lemmas required to prove Theorem 1.

\begin{lemma}\label{lemma1_1}
Let $a(n)$ be a totally multiplicative function such that $|a(p)|\leq 1$, for all $p$. Then for $x\geq 2$ and $\sigma \geq 1$, we have
$$\sum_{1< n \leq x}\frac{a(n)\Lambda(n)}{n^s\log n}=\sum_{p\leq x}\log\bigg(1-\frac{a(p)}{p^s}\bigg)^{-1}+O\bigg(\frac{1}{\sqrt x\log x}\bigg),$$
and
$$
\sum_{1\leq n \leq x}\frac{\Lambda(n)a(n)}{n^s}=\sum_{1\leq n \leq x}\frac{a(p)\log p}{p^s-a(p)}+O\bigg(\frac{1}{\sqrt{x}}\bigg).
$$
If $|a(p)|=1$ then, 
$$
\sum_{1\leq n \leq x}\frac{\Lambda(n)a(n)}{n^s}=\sum_{1\leq n \leq x}\frac{\log p}{\bar a(p) p^s-1}+O\bigg(\frac{1}{\sqrt{x}}\bigg).
$$

\end{lemma}
\begin{proof}
This follows from Lemma 5 of \cite{paper}.
\end{proof}
\begin{lemma}\label{lemma2_1}
Let $\chi$ be a primitive character mod $q$. Then, for $\Re s\geq 1$, we have,
$$\Re\sum_{1<n\leq x}\frac{\Lambda(n)\chi(n)}{n^s\log n}\leq \log\log x +C_0-\log \frac{q}{\phi(q)}+O\bigg(\frac{1}{\log x}\bigg).$$
\begin{proof}
Theorem 2.7 of \cite{classical} states that 
\begin{equation}\label{eq31}
  \sum_{{1< n\leq x}}\frac{\Lambda(n)}{n\log n}=\log\log x
+C_0+O\bigg(\frac{1}{\log x}\bigg).  
\end{equation}

By triangle inequality and \eqref{eq31} we get,
$$\Re\sum_{1< n\leq x}\frac{\Lambda(n)\chi(n)}{n^s\log n}\leq 
\sum_{{1< n\leq x\atop (q,n)=1}}\frac{\Lambda(n)}{n\log n}=\log\log x
+C_0-\log \frac{q}{\phi(q)}+O\bigg(\frac{1}{\log x}\bigg).$$
\end{proof}
\end{lemma}

\begin{lemma}\label{lemma3_1}
Assume GRH. Let $s=\sigma +\di t$ with $\sigma\geq 1$ and  let $4\leq T\leq t\leq 2T.$
Then we have
\begin{equation}
  -\frac{L'}{L}(s,\chi)=\sum_{1\leq n\leq \log^2 T} \frac{\Lambda(n)\chi(n)}{n^s}+O(\log^{2-2\sigma}T) . \label{eq1_1}  
\end{equation}

\end{lemma}
\begin{proof}
For given $x,y\geq 1$, define 
$$
 w(u)=
\begin{cases}
	1  &\qquad \hbox{if}\qquad 1\leq u\leq x,\\
	1  -\frac{\log(u/x)}{\log y} &\qquad  \hbox{if}\qquad  x<u\leq xy,\\
	 0 &\qquad  \hbox{if}\qquad  u>xy.				
\end{cases}
$$

Using Perron's formula, for $x\geq 2$ and $y\geq 2$, we have $$\sum_{n\leq xy}w(n)\frac{\Lambda(n)\chi(n)}{n^s}=\frac{-1}{2\pi\di \log y}\int_{\sigma_0-\di \infty}^{\sigma_0+\di \infty}\frac{L'}{L}(s+w)\frac{(xy)^w-x^w}{w^2}dw,$$
where $\sigma_0=\Re w.$

On pulling the contour to the left we see that, if $\chi(-1)=-1$, this is
\begin{equation}\label{logderodd}
    \frac{L'}{L}(s,\chi)=-\sum_{n\leq xy} w(n)\frac{\Lambda(n)\chi(n)}{n^s}-\sum_{\rho}\frac{(xy)^{\rho-s}-x^{\rho-s}}{(\rho-s)^2\log y}-\sum_{k=1}^{\infty}\frac{(xy)^{-2k+1-s}-x^{-2k+1-s}}{(2k-1+s)^2\log y},
\end{equation}

where $\rho=\beta+\di\gamma$ is a non trivial zero of $L(s,\chi)$ .
When $\chi(-1)=1$, we get
\begin{equation}\label{logdereven}
  \frac{L'}{L}(s,\chi)=-\sum_{n\leq xy} w(n)\frac{\Lambda(n)\chi(n)}{n^s}-\sum_{\rho}\frac{(xy)^{\rho-s}-x^{\rho-s}}{(\rho-s)^2\log y}-\sum_{k=0}^{\infty}\frac{(xy)^{-2k-s}-x^{-2k-s}}{(2k+s)^2\log y}.
  \end{equation}

We choose $y=2$ and $x=\log^2 T$ throughout the rest of the proof.
The last term on the right-hand side in \eqref{logderodd} and \eqref{logdereven} is bounded irrespective of $T$, thus it is $O(1).$\\
For simplifying the second term on the right-hand side,
we need GRH. Substituting $y=2$ and then considering the absolute values, we obtain $$\sum_{\rho}\frac{(xy)^{\rho-s}-x^{\rho-s}}{(\rho-s)^2\log y}\ll\sum_{\rho}\frac{x^{1/2-\sigma}}{|\beta-\sigma|^2+|\gamma-t|^2}.$$
GRH implies that when $L(\rho,\chi)=0$ and $\rho$ is not real then for a primitive character $\chi$ we have  $\Re\rho=\beta=1/2$. Thus for any $\sigma$ and any $\rho$ , $|\beta-\sigma|\geq 1/2.$ Also, for any $T$ and a fixed $q$, the number of roots of a Dirichlet $L$-function having ordinates between $[T,T+1]$ is  $\ll \log T.$
Thus, for $4\leq T\leq t\leq 2T$, we get $$\sum_{\gamma}\frac{1}{|\beta-\sigma|^2+|\gamma-t|^2}\ll \sum_{\gamma}\frac{1}{1+|\gamma-t|^2}\ll \log t\ll \log T.$$
Combining the above estimates, we get
\begin{equation}
\sum_{\rho}\frac{x^{1/2-\sigma}}{|\rho-s|^2}\ll x^{1/2-\sigma}{\log T}\ll \log^{2-2\sigma}T.
\end{equation}
Consider the first term on the right-hand side of $\eqref{logderodd}$. Since $y=2$, we get
\begin{equation}
\begin{split}
	\sum_{n\leq xy} w(n)\frac{\Lambda(n)\chi(n)}{n^s}&=\sum_{n\leq x}\frac{\Lambda(n)\chi(n)}{n^s}+ \sum_{x < n\leq 2x}w(n)\frac{\Lambda(n)\chi(n)}{n^s}\\
	&=\sum_{n\leq x}\frac{\Lambda(n)\chi(n)}{n^s}+ O\big(x^{-\sigma}\sum_{x< n\leq 2x}\Lambda(n)\big)\\
	&=\sum_{n\leq x}\frac{\Lambda(n)\chi(n)}{n^s}+ O(x^{1-\sigma})\\
	&=\sum_{n\leq \log^2 T}\frac{\Lambda(n)\chi(n)}{n^s}+ O(\log^{2-2\sigma}T).\\
\end{split}
\end{equation}
The last equality we get because $x=\log^2 T$.

Combining all of the above estimates, for $4\leq T\leq t\leq 2T$ and $\sigma \geq 1$, we have
$$-\frac{L'}{L}(s,\chi)=\sum_{1\leq n\leq \log^2 T} \frac{\Lambda(n)\chi(n)}{n^s}+O\big(\log^{2-2\sigma} T\big).$$
\end{proof}
\begin{lemma}\label{lemma4_1}
If $s=\sigma+\di t$ is a point such that $\sigma\geq 1$
and $4\leq T\leq t\leq 2T$, then assuming GRH we have,
$$\log L(s,\chi)= \sum_{1< n\leq \log^2 T} \frac{\Lambda(n)\chi(n)}{n^s \log n}+O\bigg({\frac{1}{\log \log T}}\bigg).$$
\end{lemma}
\begin{proof}
Since $$\log L(s,\chi)= -\int_{\sigma}^{\infty}\frac{L'}{L}(\alpha+\di t,\chi) d\alpha,$$
by substituting \eqref{eq1_1} in the above equation we get the desired result.
\end{proof}
\begin{center}
    We now proceed to prove Theorem 1.
\end{center}
\section{\textbf{Proof of Theorem 1}}
\begin{proof}
Substituting $\rho'$ for $s$ in Lemma \ref{lemma4_1} and using Lemma \ref{lemma2_1}, we get
\begin{align*}
    \log |L(\rho',\chi)|&= \Re\sum_{1< n\leq \log^2 T} \frac{\Lambda(n)\chi(n)}{n^{\rho'} \log
    n}+O\bigg({\frac{1}{\log \log T}}\bigg)\\
    &\leq\log\log (\log^2 T)+C_0-\log \frac{q}{\phi(q)} +O\bigg({\frac{1}{\log \log T}}\bigg),
\end{align*}   
where $4\leq T\leq \gamma' \leq 2T.$

Thus, $$\frac{|L(\rho',\chi)|}{\log\log T}\leq 2e^{C_0}\frac{\phi(q)}{q}\bigg(1+O\bigg(\frac{1}{\log\log T}\bigg)\bigg).$$
Since $4\leq T\leq \gamma'\leq 2T,$
$$\limsup\limits_{\gamma' \rightarrow\infty}\frac{|L(\rho',\chi)|}{\log\log \gamma'}\leq 2e^{C_0}\frac{\phi(q)}{q} .$$
\end{proof}
Remark: In the proof of Theorem 1 we did not use any special property of critical points of the zeta function. One can also recover
the following classical bound.\\
Let $s=\sigma +\di t$ be a point with $\sigma\geq1$, then
$$\limsup\limits_{t \rightarrow\infty}\frac{|L(s,\chi)|}{\log\log t}\leq 2e^{C_0}\frac{\phi(q)}{q}.$$

\section{\textbf{Lemmas related to Theorem 2}}
In this section we prove lemmas required to prove Theorem 2.

Let $\delta$ be such that $1/2<\delta<1$. Define $\epsilon=2\delta-1$.
 Let $x$ be a sufficiently large number such that for the given $\epsilon$, we have $x>q^{1/\epsilon}$. Let $m$ be some positive number greater than 1. We shall determine the exact value of $m$ later, while $q$ remains fixed throughout. For $x>q^{1/\epsilon}$, define a function $$W_x(s)=\sum_{n\leq x}\frac{\Lambda(n)b(n)}{n^s},$$ where $b(n)$ is a completely multiplicative function such that  \\
$$
b(p)=\left.
\begin{cases}
1 &\text {if p $\mid$ q},\\
\bar \chi(p) &\text {if $p\leq x^\epsilon$ and $p \nmid q$},\\
-1& \text{if $x^\epsilon<p\leq mx^\epsilon$},\\
1 & \text{if $mx^\epsilon<p\leq x^\delta$},\\
-1& \text{if $x^\delta<p$}.
\end{cases}
\right.
$$
For an appropriately large choice of $x$, define $$V_x(s)=\sum_{n\leq x}\frac{\Lambda(n)}{n^s}.$$
Using $V_x(s)$ and $W_x(s)$ as auxiliary functions, we will construct roots (zeros) of the logarithmic derivative of the Riemann zeta function to the right of the vertical line $\Re s=1$.
Without assuming RH and by considering just the standard zero-free region and Perron's formula with $|s-1|<1/10$, for some $\tilde{c}>0$, we obtain
\begin{equation}\label{eq1_2}
 V_x(s)=\frac{x^{1-s}}{1-s}-\frac{\zeta'}{\zeta}(s)+O\big(x^{1-\sigma}\exp(\tilde{-c}\sqrt{\log x})\big). 
\end{equation}
We use the above equation to determine the behaviour of $W_x(s)$ near $s=1$. 
 
 Define 
 \begin{equation}\label{eqs1}
   S_1=-\frac{L'}{L}(1,\bar\chi)+\sum_{p \mid q}(\log p )/(p-1).  
 \end{equation}
 
Define $m$ such that $2\log m=\Re S_1+|S_1|+1.$ 

\begin{lemma}\label{lemma1_2}
If $x^\epsilon>q$, $\Re s\geq 1$ and $$|s-1|\leq\frac{6|S_1|+3}{(-2\delta^2+1+\epsilon^2)\log^2 x},$$ then 
$$W_x(s)=\frac{(1-s)}{2}(2\delta^2-1-\epsilon^2)\log^2 x+S_1-2\log m+O\bigg(\frac{1}{\log x}\bigg).$$ From above, we see that $W_x(s) $ has a root at $$s=1+\frac{2(S_1-2\log m)}{(2\delta^2-1-\epsilon^2)\log^2 x}+O\bigg(\frac{1}{\log^3 x}\bigg).$$
\end{lemma}
\begin{proof}
Our aim here is to represent $W_x(s)$ in terms of $V_x(s)$, and to show that $W_x(s)$ has a root near the point $s=1$ and to the right of the vertical line $\Re s=1$. \\
First, we split $W_x(s)$ into four different parts, as follows:
\begin{equation}\label{eqfirst}
    \begin{split}
         W_x(s)=\sum_{1\leq n\leq x^\epsilon}\frac{\Lambda(n)b(n)}{n^s}+\sum_{x^\epsilon< n\leq mx^\epsilon}\frac{\Lambda(n)b(n)}{n^s}+\sum_{mx^\epsilon< n\leq x^\delta}\frac{\Lambda(n)b(n)}{n^s}+\sum_{x^\delta< n\leq x}\frac{\Lambda(n)b(n)}{n^s}.
         \end{split}
  \end{equation}
Let us represent the four sums by $T_1, T_2, T_3,T_4$ respectively.
We know the following from Mertens' formulae and the prime number theorem. For $y\geq 1$ and for some constant $r$, we have
\begin{align}
    \sum_{p\leq y}\frac{\log p}{p}&=\log y+r+O(\exp(-c\sqrt{\log y}))\label{eqformula1},\\
    \sum_{p\leq y}\frac{\log^2 p}{p}&=\frac{1}{2}\log^2 y + O(\log y).\label{eqformula2}
\end{align}
Also, for $a>1$ and $1<y_1<y_2$, we have
\begin{align}\label{eqap}
  \sum_{y_1\leq p\leq y_2}\frac{\log p}{p^a}=\int_{y_1}^{y_2}\frac{d\theta(x)}{x^a}=O(y_1^{1-a}),  
\end{align}

where $\theta(x)=\sum_{p\leq x}\log p.$\\

First we estimate $T_1$. 

Using Perron's formula, for $\Re s\geq 1$ and $x\geq 1$, we have the following result.
\begin{equation}\label{eqn57}
\begin{split}
     -\frac{L'}{L}(s,\bar\chi)&=\sum_{n\leq x^\epsilon}\frac{\Lambda(n)\bar\chi(n)}{n^s}+O\big(x^{(1-\sigma)\epsilon}\exp(-\tilde c\sqrt{\epsilon\log x})\big)\\
    &=\sum_{n\leq x^\epsilon}\frac{\Lambda(n)b(n)}{n^s}+\sum_{n\leq x^\epsilon}\frac{\Lambda(n)\bar\chi(n)}{n^s}- \sum_{n\leq x^\epsilon}\frac{\Lambda(n)b(n)}{n^s}+O\big(x^{(1-\sigma)\epsilon}\exp(-\tilde c\sqrt{\epsilon\log x})\big)\\
    &=T_1-\sum_{n\leq x^\epsilon \atop (n,q)>1}\frac{\Lambda(n)}{n^s}+O\big(1/\log x\big).
\end{split}
 \end{equation}
We see that,
\begin{equation}\label{eqn58}
    \begin{split}
  \sum_{n\leq x^\epsilon \atop(q,n)>1}\frac{\Lambda(n)}{n^s}-\sum_{n\leq x^\epsilon \atop(q,n)>1}\frac{\Lambda(n)}{n}&=\sum_{p^k\leq x^\epsilon \atop p|q}\frac{\log p}{p^k}\bigg(e^{k(1-s)\log p}-1\bigg)      \\
  &=\sum_{p^k\leq x^\epsilon \atop p|q}\frac{\log p}{p^k}\bigg(\sum_{j=1}^{\infty}\frac{(k(1-s)\log p)^j}{j!}\bigg)\\  
  &\ll \frac{1}{\log x}\sum_{p\leq x^{\epsilon/k} \atop p|q}\frac{\log p}{p^k}\ll \frac{1}{\log x}.
\end{split}
\end{equation}

Also,
 \begin{equation}\label{eqn59}
     \begin{split}
         \sum_{n\leq x^\epsilon\atop (q,n)>1}\frac{\Lambda(n)}{n}&=\sum_{n=1\atop (q,n)>1}^{\infty}\frac{\Lambda(n)}{n}-\sum_{n> x^\epsilon\atop (q,n)>1}\frac{\Lambda(n)}{n}=\sum_{p \mid q}(\log p )/(p-1)+O\bigg(1/\log x\bigg).
     \end{split}
 \end{equation}
By plugging in \eqref{eqn58} and \eqref{eqn59} in \eqref{eqn57}, we have
\begin{equation*}
-\frac{L'}{L}(s,\bar\chi)=T_1- \sum_{p \mid q}(\log p )/(p-1)+O\big(1/\log x\big).\\
\end{equation*}
Let $D_x(s)=\sum_{n\leq x^\epsilon}{\Lambda(n)\bar\chi(n)}{n^{-s}}.$
Because $D'_x(s)=\sum_{n\leq x^\epsilon}-{(\log n)\Lambda(n)\bar\chi(n)}{n^{-s}}\ll1,$
we have $$D_x(s)=D_x(1)+O(1/\log x).$$
So we get
$$-\frac{L'}{L}(s,\bar\chi)=D_x(s)+ O\big(x^{(1-\sigma)\epsilon}\exp(-\tilde c\sqrt{\epsilon\log x})\big),$$
and
$$-\frac{L'}{L}(1,\bar\chi)=D_x(1)+ O\big(\exp(-\tilde c\sqrt{\log x})\big),$$

and we conclude
\begin{equation}\label{newbound}
  -\frac{L'}{L}(s,\bar\chi)=-\frac{L'}{L}(1,\bar\chi)+O\big(1/\log x\big).  
\end{equation}

Thus,
 \begin{equation}
     \begin{split}
          T_1 &= -\frac{L'}{L}(s,\bar\chi)+\sum_{p \mid q}(\log p )/(p-1)+O\big(1/\log x\big).\\
          \end{split}
\end{equation}

Recall that $$S_1=-\frac{L'}{L}(1,\bar\chi)+\sum_{p \mid q}(\log p )/(p-1).$$ 
So we have
\begin{equation}\label{eq3_2}
T_1=S_1 +O\big(1/\log x).
\end{equation}

We estimate $T_2$ as follows.
\begin{align*}
T_2=\sum_{x^\epsilon< n\leq mx^\epsilon}\frac{\Lambda(n)b(n)}{n^s}&=-\sum_{x^\epsilon< p\leq mx^\epsilon}\frac{\log p}{p^{s}}+O\bigg(\sum_{k=2 }^{\log(mx^\epsilon)/\log 2} \qquad \sum_{p=x^{\epsilon/k}}^{(mx^\epsilon)^{1/k}}\frac{\log p}{p^{k}}\bigg).\\ \end{align*}
 
By \eqref{eqap} we see that $\sum_{p=x^{\epsilon/k}}^{(mx^\epsilon)^{1/k}}\frac{\log p}{p^{k}}= O(x^{\epsilon (1-k)/k}).$
Thus,$$\sum_{k=2 }^{\log(mx^\epsilon)/\log 2} \qquad \sum_{p=x^{\epsilon/k}}^{(mx^\epsilon)^{1/k}}\frac{\log p}{p^{k}}= \sum_{k=2 }^{\log(mx^\epsilon)/\log 2}O(x^{\epsilon (1-k)/k})=O\bigg(\frac{1}{\log x}\bigg).$$
Thus, $$T_2=-\sum_{x^\epsilon< p\leq mx^\epsilon}\frac{\log p}{p^{s}}+O\bigg(\frac{1}{\log x}\bigg).$$

We also see that,
\begin{equation}\label{eq5.9new}
\begin{split}
    \frac{\log p}{p^s}-\frac{\log p}{p}=(e^{\log p(1-s)}-1)\ll\frac{\log p}{p}\frac{\log p}{\log^2 x}.\\ 
    \end{split}
\end{equation}
Using \eqref{eqformula2} and \eqref{eq5.9new}, we have
\begin{align*}
    \sum_{x^\epsilon< p\leq mx^\epsilon}\bigg(\frac{\log p}{p^s}-\frac{\log p}{p}\bigg)&\ll \sum_{x^\epsilon< p\leq mx^\epsilon}\frac{\log^2 p}{p \log ^2 x}\ll\frac{1}{\log^2 x}\sum_{x^\epsilon< p\leq mx^\epsilon}\frac{\log^2 p}{p} \ll \frac{1}{\log x}.
\end{align*}
 Thus, from \eqref{eqformula1}we have 
 \begin{equation}\label{eq4_2}
  T_2= -\sum_{x^\epsilon< p\leq mx^\epsilon}\frac{\log p}{p}+O(\log^{-1}x)=-\log m +O(\log^{-1}x).   
 \end{equation}

We estimate $T_3$ and $T_4$ as shown below.
\begin{equation}\label{eq5_2}
    \bigg|T_3-\sum_{mx^\epsilon <n\leq x^\delta}\frac{\Lambda(n)}{n^s}\bigg|\ll \sum_{k=2 }^{\delta\log x/\log 2}\sum_{p=(mx^\epsilon)^{1/k}}^{x^{\delta/k}}\frac{\log p}{p^k}\ll \sum_{k=2 }^{\delta\log x/\log 2}\frac{1}{{(mx^{\epsilon})}^\frac{k-1}{k}}\ll(\log x)^{-1}.
\end{equation}
Thus,
\begin{equation}\label{eq5_2n}
   T_3=V_{x^\delta}(s)-V_{mx^\epsilon}(s)+O(1/\log x).
\end{equation}

\begin{equation}\label{eq6_2}
\bigg|T_4+\sum_{x^\delta <n\leq x}\frac{\Lambda(n)}{n^s}\bigg|\ll \sum_{k=2 }^{\log x/\log 2}\sum_{p=x^{\delta/k}}^{x^{1/k}}\frac{\log p}{p^k}\ll \sum_{k=2 }^{\log x/\log 2}\frac{1}{x^\frac{(k-1){\delta}}{k}}\ll (\log x)^{-1}.    
\end{equation}
Thus,
\begin{equation}\label{eq6_2n}
    T_4=V_{x^\delta}(s)-V_{x}(s)+O(1/\log x).
\end{equation}

The third inequality in \eqref{eq5_2} and \eqref{eq6_2} is due to \eqref{eqap}. Using \eqref{eqfirst}, \eqref{eq3_2}, \eqref{eq4_2},\eqref{eq5_2n}, \eqref{eq6_2n}, we get 

\begin{equation}\label{eqw}
W_x(s)=2V_{x^\delta}(s)-V_x(s)-V_{mx^\epsilon}(s)+S_1-\log m+O\bigg(\frac{1}{\log x}\bigg).    
\end{equation}

Now, we substitute a suitable version of \eqref{eq1_2} in \eqref{eqw} to get
\begin{equation*}\label{weqation}
\begin{split}
W_x(s)&=\frac{2x^{\delta(1-s)}}{1-s}-\frac{x^{1-s}}{1-s}-\frac{(mx^{\epsilon})^{(1-s)}}{1-s}+S_1-\log m+O\bigg(\frac{1}{\log x}\bigg)\\
&=\bigg(\frac{2e^{\delta(1-s)\log x}-e^{(1-s)\log x}-e^{(\log m+\epsilon\log x)(1-s)}}{1-s}\bigg)+S_1-\log m+O(1/\log x)\\
&= \frac{(1-s)}{2}\big[(2\delta^2-1-\epsilon^2)\log^2 x\big]-\frac{1-s}{2}(\log^2m+2\epsilon \log m\log x )+S_1-2\log m+O\bigg(\frac{1}{\log x}\bigg).
\end{split}
\end{equation*}
Now since, $$|s-1|\leq\frac{6|S_1|+3}{(-2\delta^2+1+\epsilon^2)\log^2 x},$$
we have
$$W_x(s)=\frac{(1-s)}{2}(2\delta^2-1-\epsilon^2)\log^2 x+S_1-2\log m+O\bigg(\frac{1}{\log x}\bigg).
$$
Thus, $W_x(s) $ has a root at $$s=1+\frac{2(S_1-2\log m)}{(2\delta^2-1-\epsilon^2)\log^2 x}+O\bigg(\frac{1}{\log^3 x}\bigg).$$

Note that $(2\delta^2-1-\epsilon^2)<0$ for any positive choice of $\epsilon$ and $\delta$ that satisfy the conditions mentioned at the beginning of this section.
For the given choice of $m$, we can see that $(\Re S_1-2\log m)<0.$ Since $2(\Re S_1-2\log m)/(2\delta^2-1-\epsilon^2)>0$, we can say that the root has the real part greater than 1.
\end{proof}
\begin{center}
We now proceed to prove Theorem 2.
\end{center}

\section{\textbf{Proof of Theorem 2}}
\begin{proof}
Define circles $\mathscr{C}_0$ and $\mathscr{C}_1$ as
$$\mathscr{C}_0= \bigg\{s=1+\frac{2(-S_1+2\log m)}{(-2\delta^2+1+\epsilon^2)\log^2 x} +\frac{(-\Re S_1+2\log m)}{(-2\delta^2+1+\epsilon^2)\log^2 x}e^{\di\phi}:0\leq \phi\leq 2\pi\bigg\},$$
$$\mathscr{C}_1= \bigg\{s=1+\frac{2(-S_1+2\log m)}{(-2\delta^2+1+\epsilon^2)\log^2 x} +\frac{(-\Re S_1+2\log m)}{2(-2\delta^2+1+\epsilon^2)\log^2 x}e^{\di\phi}:0\leq \phi\leq 2\pi\bigg\}.$$
Let $\|y\|$ represent the distance of a real number $y$ from the nearest integer. For a given $x$, let $\tau$ be a number such that, 
\begin{equation}\label{eq514}
\|\tau\log p/2\pi-\arg{ \bar b(p)/2\pi}\|\leq1/\log^2 x     
\end{equation}
for all $p\leq x$. Here we consider $\arg \bar b(p) \mod 2\pi$.
For such $\tau$ and for all $p\leq x$ using equation (35) from \cite{paper}, we get
\begin{equation}\label{eq8_2}
|p^{-\di\tau}-b(p)|\ll 1/\log^2 x.
\end{equation}

For a $\tau$ satisfying condition \eqref{eq8_2},
and for $s \in \mathscr{C}_1$ and for all $p\leq x$, we have
\begin{equation}
    \begin{split}
    |V_x(s+\di\tau)-W_x(s)|&\leq\sum_{ n \leq x}\bigg|\frac{\Lambda(n)}{n^{s+\di\tau}}-\frac{\Lambda(n)b(n)}{n^s}\bigg|\ll\sum_{p^k \leq x}\frac{k\log p}{p^k}|p^{-\di\tau}-b(p)|\ll\frac{1}{\log x}.
    \end{split}
  \end{equation}
Thus,
\begin{equation}\label{eq13_1}
    V_x(s+\di\tau)=W_x(s)+O\bigg(\frac{1}{\log x}\bigg).
\end{equation}

For $x>q^{1/\epsilon}$, define $$M_x(s)=\sum_{1<n\leq x}\frac{\Lambda(n)c(n)}{n^s\log n},$$
where $c(n)$ is a completely multiplicative function such that 
$$
c(p)=
\begin{cases}
 1 &\text{if\qquad} p \mid q,\\
\chi(p)b(p) &\text{if\qquad}p\nmid q.\\
 \end{cases}
$$
We define $c(p)$ as above because for a $\tau$ satisfying \eqref{eq8_2} and $p\nmid q$, we want $\chi(p)p^{-\di \tau}\approx c(p)$. But we know that $\chi(p)p^{-\di \tau}\approx \chi(p)b(p)$.Thus, $c(p)=\chi(p)b(p)$ for $p\nmid q$.

By Lemma \ref{lemma1_1}, when $s\in \mathscr{C}_1$, we get
\begin{equation}
\begin{split}
\sum_{1<n\leq x \atop (q,n)>1}\frac{\Lambda(n)c(n)}{n^s\log n}&=\sum_{p \mid q }\log\bigg(1-\frac{c(p)}{p^s}\bigg)^{-1}+O\bigg(\frac{1}{\sqrt{x}\log x}\bigg)\\
&=-\sum_{p \mid q }\log\bigg(1-\frac{1}{p^s}\bigg)+O\bigg(\frac{1}{\sqrt{x}\log x}\bigg).  \\
\end{split}
\end{equation}

For $s\in \mathscr{C}_1$, we have $p^{1-s}=1+ O\big(1/\log^2 x\big)$ and $\frac{p^s-1}{p-1}=1+O\big(1/\log^2 x\big).$
So we have
\begin{align*}
\log\bigg(1-\frac{1}{p^s}\bigg)-\log\bigg(1-\frac{1}{p}\bigg)=\log\bigg(\frac{p^{1-s}(p^s-1)}{(p-1)}\bigg)\ll\frac{1}{\log^2 x}.
\end{align*}
Consequently, we obtain
\begin{equation}\label{eqlog}
    \sum_{1<n\leq x \atop (q,n)>1}\frac{\Lambda(n)c(n)}{n^s\log n}=-\sum_{p \mid q }\log\bigg(1-\frac{1}{p}\bigg)+O\bigg(\frac{1}{\log^2 x}\bigg).
\end{equation}

Thus,
\begin{equation}
     \begin{split}
         &\sum_{1<n\leq x}\frac{\Lambda(n)\chi(n)}{n^{s+\di\tau}\log n}-M_x(s)-\log\frac{\phi(q)}{q}=\sum_{1<n\leq x}\frac{\Lambda(n)\chi(n)}{n^{s+\di\tau}\log n}-\sum_{1<n\leq x}\frac{\Lambda(n)c(n)}{n^s\log n}+\sum_{p\mid q}\sum_{m=1}^{\infty}\frac{1}{mp^m}.\\
\end{split}
\end{equation}

Using \eqref{eqlog}, for any $\tau$ satisfying \eqref{eq8_2} and $s\in \mathscr{C}_1$, we see that the right-hand side of the above equation is
\begin{equation}\label{boundL}
\begin{split}
    &\leq \sum_{1<n\leq x \atop (q,n)=1}\frac{\Lambda(n)}{n\log n}\bigg|\frac{\chi(n)}{n^{\di\tau}}-c(n)\bigg|+O\bigg(\frac{1}{\log x}\bigg)
     \leq \sum_{1<p^k\leq x \atop p\nmid q}\frac{1}{kp^k}\bigg|\frac{\chi(p^k)}{p^{k\di\tau}}-c(p^k)\bigg|+O\bigg(\frac{1}{\log x}\bigg)\\
     &\leq \sum_{1<p^k\leq x\atop p\nmid q }\frac{1}{p^k}\bigg|\frac{\chi(p)}{p^{\di\tau}}-c(p)\bigg|+O\bigg(\frac{1}{\log x}\bigg).   
\end{split}
 \end{equation}

For $p\leq x$ and $p\nmid q$, using \eqref{eq8_2}, we have 
\begin{equation*}
  |p^{-\di \tau}-b(p)|=|\chi(p)p^{-\di \tau}-\chi(p)b(p)|=\bigg|\frac{\chi(p)}{p^{\di\tau}}-c(p)\bigg|\leq \frac{1}{\log^2 x}. 
  \end{equation*}

So we get
\begin{equation}\label{est1}
     \sum_{1<p^k\leq x \atop p\nmid q}\frac{1}{p^k}\bigg|\frac{\chi(p)}{p^{\di\tau}}-c(p)\bigg|\ll \sum_{1<p\leq x^{1/k} \atop p\nmid q}\frac{1}{p^k}\bigg|\frac{\chi(p)}{p^{\di\tau}}-c(p)\bigg|\ll\frac{1}{\log x}.
\end{equation}

Combining \eqref{boundL} and \eqref{est1}, we get
\begin{equation}\label{eq10_2}
  \sum_{1<n\leq x}\frac{\Lambda(n)\chi(n)}{n^{s+\di\tau}\log n}=M_x(s)+ \log \frac{\phi(q)}{q} +O\bigg(\frac{1}{\log x}\bigg).  
\end{equation}
Also,
\begin{equation}\label{Lderivative}
 M'_x(s)\ll \sum_{1\leq n \leq x}\frac{\Lambda(n)c(n)}{n}\ll \log x.   
\end{equation}
Thus, from \eqref{Lderivative} and using Liouville's theorem, we get 
\begin{equation}\label{eq11_2}
  M_x(s)=M_x(1)+O\bigg(\frac{1}{\log x}\bigg).  
\end{equation}

We have  $$M_x(1)=\sum_{1<n\leq x^\epsilon}\frac{\Lambda(n)c(n)}{n\log n}+\sum_{x^\epsilon<n\leq mx^\epsilon}\frac{\Lambda(n)c(n)}{n\log n}+\sum_{mx^\epsilon<n\leq x^\delta}\frac{\Lambda(n)c(n)}{n\log n}+\sum_{x^\delta<n\leq x}\frac{\Lambda(n)c(n)}{n\log n}.$$

We estimate the second, third, and fourth terms on the right-hand side of the above equation as follows:
\begin{equation}\label{eq15_2}
  \sum_{x^\epsilon<n\leq mx^\epsilon}\frac{\Lambda(n)c(n)}{n\log n}=-\sum_{x^\epsilon <p\leq mx^\epsilon}\frac{\chi(p)}{p}+O\bigg(\sum_{k=2 }^{\log (mx)^\epsilon/\log 2}\qquad \sum_{p=x^{\epsilon/k}}^{(mx^\epsilon)^{1/k}}\frac{1}{p^k}\bigg)\ll \frac{1}{\log x},  
\end{equation}

and
\begin{equation}\label{eq16_2}
    \sum_{mx^\epsilon<n\leq x^\delta}\frac{\Lambda(n)c(n)}{n\log n}=\sum_{mx^\epsilon<p\leq x^\delta}\frac{\chi(p)}{p}+O\bigg(\sum_{k=2 }^{\delta \log x/\log 2}\qquad \sum_{p=(mx^\epsilon)^{1/k}}^{x^{\delta/k}}\frac{1}{p^k}\bigg)\ll \frac{1}{\log x},
\end{equation}    
and
\begin{equation}\label{eq17_2}
\sum_{x^\delta<n\leq x}\frac{\Lambda(n)c(n)}{n\log n}=-\sum_{x^\delta<p\leq x}\frac{\chi(p)}{p}+O\bigg(\sum_{k=2 }^{\log x/\log 2}\qquad \sum_{p=x^{\delta/k}}^{x^{1/k}}\frac{1}{p^k}\bigg)\ll \frac{1}{\log x}.
\end{equation}
The first term on the right-hand side of the above sums is approximated using partial summation.
Following is the estimate for $M_x(1),$ which we get from 
\eqref{eq15_2}, \eqref{eq16_2} and \eqref{eq17_2}.

\begin{equation}\label{eq12_2}
\begin{split}
  M_x(1)&=\sum_{1<n\leq x}\frac{\Lambda(n)c(n)}{n\log n}=\sum_{1<n\leq x^\epsilon}\frac{\Lambda(n)}{n\log n}+O\bigg(\frac{1}{\log x}\bigg)\\
&=\log\log x^\epsilon +C_0+O\bigg(\frac{1}{\log x}\bigg) \\
&=\log\log x +\log \epsilon +C_0 +O\bigg(\frac{1}{\log x}\bigg).
\end{split}    
\end{equation}
The second line in the above equation follows from Theorem 2.7 of \cite{classical}.

Combining \eqref{eq10_2}, \eqref{eq11_2} and \eqref{eq12_2}, we obtain
\begin{equation}\label{614}
  \Re \sum_{1<n\leq x}\frac{\Lambda(n)\chi(n)}{n^{s+\di\tau}\log n}\geq \log\log x +C_0 +\log \frac{\phi(q)}{q} +\log \epsilon+O\bigg(\frac{1}{\log x}\bigg).  
\end{equation}
For $x \geq 1$ and $\Re s>1$, define $$T_x(s)=\sum_{p>x}\frac{\chi(p)}{p^s},$$
and $$Q_x(s)=\sum_{p>x}\frac{1}{p^s}.$$ Let $Q'_x(s)$ be the term by term derivative of $Q_x(s)$ with respect to $s$.\\

For $\Re s>1 $, we have
$$\log L(s,\chi)= \sum_{1<n\leq x}\frac{\Lambda(n)\chi(n)}{n^{s}\log n}+T_x(s)+O\bigg(\frac{1}{\log x}\bigg),$$ and
$$-\frac{\zeta'}{\zeta}(s)=\sum_{1\leq n\leq x}\frac{\Lambda(n)}{n^s}-Q'_x(s)+O\bigg(\frac{1}{\log x}\bigg).$$
If $s\in \mathscr{C}_1$ and $\tau$ satisfying \eqref{eq8_2} is such that 
\begin{equation}\label{eqT}
    |T_x(s+\di\tau)|\leq 1/\log x,
\end{equation}
and
\begin{equation}\label{eqQ}
    |Q'_x(s+\di\tau)|\leq 1/\log x
\end{equation}
then we get
\begin{equation}\label{eq618n}
\Re \log L(s+\di\tau,\chi) \geq \log \log x+C_0+\log \frac{\phi(q)}{q}+\log \epsilon + O\bigg(\frac{1}{\log x}\bigg),
\end{equation}
and
\begin{equation}\label{rouche}
-\frac{\zeta'}{\zeta}(s+\di\tau)=W_x(s)+O\bigg(\frac{1}{\log x}\bigg).   \end{equation}

If $W_x(s)\asymp 1$ on the boundary of $\mathscr{C}_1$, then
\begin{equation}\label{rouche1}
-\frac{\zeta'}{\zeta}(s+\di\tau)-W_x(s)\ll\frac{1}{\log x}\ll W_x(s).   
\end{equation}
From \eqref{rouche1}, we can see that if $W_x(s)$ has a root for $s \in \mathscr{C}_1$ and $W_x(s)\asymp 1$ on the boundary of $\mathscr{C}_1$, then by Rouche's theorem, $\zeta'(s)$ will have a zero $\rho'$ in $\mathscr{C}_1+\di \tau $
and from \eqref{eq618n} we deduce that
\be\label{eq5_17}
 {|L(\rho',\chi)|}\geq  \epsilon e^{C_0} \frac{\phi(q)}{q}e^{O(1/\log x)}{\log x}. 
\ee

Finally, to prove Theorem 2, we want the following to be true:
\begin{itemize}
  \item $W_x(s)$ to have a root inside $\mathscr{C}_1,$ 
  \item $W_x(s)\asymp 1$ on $\mathscr{C}_1,$ 
  \item Existence of arbitrarily large $\tau$ satisfying \eqref{eqT} and \eqref{eqQ}.
\end{itemize}
The fulfillment of the first two conditions above depends on making the correct choice for $m$. Recall that 
\begin{equation}\label{mvalue}
  2\log m=\Re S_1+|S_1|+1.  
\end{equation}

\begin{lemma}
$W_x(s)$ has a root inside $\mathscr{C}_1$. 
\end{lemma}
\begin{proof}
To prove this, we need to check that the distance of the root of $W_x(s)$ from the center of $\mathscr{C}_1$ is less than its radius.

From Lemma \ref{lemma1_2} we know that $W_x(s)$ has a root at 
$$s=1+\frac{2(S_1-2\log m)}{(2\delta^2-1-\epsilon^2)\log^2 x}+O\bigg(\frac{1}{\log^3 x}\bigg).$$
From the definition of $\mathscr{C}_1$, we can see that the distance between its zero and the center of $\mathscr{C}_1$ is $O(1/\log^3 x)$, whereas the radius of $\mathscr{C}_1$ is $$\frac{(-\Re S_1+2\log m)}{(-2\delta^2+1+\epsilon^2)\log^2 x}.$$
Thus proving the assertion.
\end{proof}

\begin{lemma}
    We have $W_x(s)\asymp 1$ on the boundary of $\mathscr{C}_1.$
\end{lemma}
\begin{proof}
To check how $W_x(s)$ varies on $\mathscr{C}_1$, we plug the equation of $\mathscr{C}_1$ in the equation of $W_x(s)$. We get,
\begin{equation*}
\begin{split}
W_x(\mathscr{C}_1)&=\frac{-1}{2}\bigg(\frac{2(-S_1+2\log m)}{(-2\delta^2+1+\epsilon^2)\log^2 x} +\frac{(-\Re S_1+2\log m)}{2(-2\delta^2+1+\epsilon^2)\log^2 x}e^{\di\phi}\bigg)(2\delta^2-1-\epsilon^2)\log^2x+S_1-2\log m\\
&+O\bigg(\frac{1}{\log x}\bigg)\\
&=\frac{1}{4}(2\log m-\Re S_1)e^{\di \phi}+O\bigg(\frac{1}{\log x}\bigg)\asymp 1.
\end{split}
\end{equation*}
    Thus, $W_x(s)\asymp 1$ on the boundary of $\mathscr{C}_1.$
\end{proof}
\begin{lemma}
    There exist arbitrarily large $\tau$ satisfying \eqref{eqT} and \eqref{eqQ}.  
\end{lemma}
\begin{proof}
If $\tau$ has a property such that 
\be\label{eq5_18}
\oint\limits_{\mathscr{C}_0}|T_x(z+\di\tau)|^2\,|\mathrm{d}z|\leq 1/\log^8 x\ee
and 
\be\label{eq5_19}
\oint\limits_{\mathscr{C}_0}|Q_x(z+\di\tau)|^2\,|\mathrm{d}z|\leq 1/\log^8 x,\ee  

then whenever $s$ lies in or on  
$\mathscr{C}_1$ by Cauchy-Schwarz inequality, we get 

\be
|T_x(s+\di\tau)|=\bigg|\frac{1}{2\pi\di}\oint\limits_{\mathscr{C}_0}\frac{T_x(z+\di\tau)}{z-s}\mathrm{d}z\bigg|\leq\log x\bigg(\oint\limits_{\mathscr{C}_0}|T_x(z+\di\tau)|^2|\mathrm{d}z|\bigg)^{1/2}\leq \frac{1}{\log x},
\ee
and
\be
|Q'_x(s+\di\tau)|=\bigg|\frac{1}{2\pi\di}\oint\limits_{\mathscr{C}_0}\frac{Q_x(z+\di\tau)}{(z-s)^2}\mathrm{d}z\bigg|\leq \log^3 x\bigg(\oint\limits_{\mathscr{C}_0}|Q_x(z+\di\tau)|^2|\mathrm{d}z|\bigg)^{1/2}\leq \frac{1}{\log x}.
\ee
To prove the existence of $\tau$ satisfying \eqref{eq5_18} and \eqref{eq5_19} and \eqref{eq8_2}, we establish some necessary setup.\\
 
We would like to show that there exists an arbitrarily large $t$ such that $\big\|{t\log p/2\pi} -\arg \bar b(p)/2 \pi \big\|\leq 1/\log^2 x$ for all primes $p\leq x$. This will imply that such $t$ satisfies $|p^{-\di t}-b(p)|\ll 1/\log^2 x$ for all $p\leq x.$ To accomplish this Montgomery and Gonek \cite{paper} defined a function $h(t)$ such that $h(t)>0$ only when 
$\big\|{t\log p/2\pi} -\arg{ \bar b(p)/2\pi} \big\|\leq 1/\log^2 x$ for all $p\leq x$ .  

Also, from Lemma 9 of Montgomery and Gonek \cite{paper} for $x=\frac{\log T}{3(\log\log T)^3}$, we have
\begin{equation}\label{eq623}
    \int_{T}^{2T}h(t)dt=(1+O(1/x))\mu^{\pi(x)}T,
\end{equation}
where $\pi(x)$ is the number of primes less than or equal to $x$.
Equation \eqref{eq623} shows that there is a positive proportion of values $t\in[T,2T]$ such that $h(t)>0$. Which implies that there is a positive proportion of values $t\in[T,2T]$ such that $|p^{-\di t}-b(p)|\ll 1/\log^2 x$ for all $p\leq x$. 

Montgomery and Gonek \cite{paper} also defined $g(t)$ with the following property.\\ 
For each $p>x$, let $b_p$ have the property that $|b_p| \leq 1/p^{\sigma}$.
Then for $\sigma>1$ we have
 \begin{equation}\label{eq624}
  \int_{T}^{2T}g(t)\bigg|\sum_{p>x}\frac{b_p}{p^{\di t}}\bigg|^2dt\ll T\mu^{\pi(x)}\bigg(\frac{1}{x}+T^{-1/2}\log\frac{\sigma}{\sigma-1}\bigg).   
 \end{equation}
The above equation follows from Lemma 10 of \cite{paper}.

Please note that here we avoid restating the definitions of $h(t)$ and $g(t)$ and only state the essential properties of these functions.

Define $h^+(t)=\max\{0,h(t)\}$.

Since $h(t)\leq h^+(t)\leq g(t)$, we obtain
$\int_{T}^{2T}h^+(t)dt=\mu^{\pi(x)}T\big(1+O(1/x)\big)$\cite{paper}.

We have the following bound.
\begin{align*}
  &\int_{T}^{2T}h^+(t)\bigg(\oint\limits_{\mathscr{C}_0}\big(|Q_x(z+\di t)|^2+|T_x(z+\di t)|^2\big)|\mathrm{d}z|\bigg)dt\\ 
  &\leq \oint\limits_{\mathscr{C}_0}\bigg(\int_{T}^{2T}g(t)\big(|Q_x(z+\di t)|^2+|T_x(z+\di t)|^2\big)dt\bigg)|dz|.\\
\end{align*}
Since $z\in \mathscr{C}_0$, $\Re z\geq 1+ c/{2\log^2 x}$, for some positive $c$. Substituting $\sigma=\Re z$ in $\eqref{eq624}$, we get 
$$
\int_{T}^{2T}g(t)|Q_x(z+\di t)|^2dt\ll T\mu^{\pi(x)}/x,
$$
and 
$$
\int_{T}^{2T}g(t)|T_x(z+\di t)|^2dt\ll T\mu^{\pi(x)}/x.
$$
Combining the above two estimates, we get
\begin{equation}
\int_{T}^{2T}h^+(t)\bigg(\oint\limits_{\mathscr{C}_0}\big(|Q_x(z+\di t)|^2+|T_x(z+\di t)|^2\big)|\mathrm{d}z|\bigg)dt\ll \oint\limits_{\mathscr{C}_0} \mu^{\pi(x)}T/x|dz|\ll\mu^{\pi(x)}Tx^{-1}\log^{-2}x.   
\end{equation}

Since the minimum of a function is at the most as large as its average, it follows that there is a $\tau$ with
$T \leq \tau \leq 2T$ and $h^+(\tau) > 0$, such that $$\oint\limits_{\mathscr{C}_0}|T_x(z+\di\tau)|^2\,|\mathrm{d}z|\leq 1/\log^8 x,$$ and $$\oint\limits_{\mathscr{C}_0}|Q_x(z+\di\tau)|^2\,|\mathrm{d}z|\leq 1/\log^8 x.$$
So there exists a $\tau$ such that it satisfies \eqref{eq8_2},\eqref{eq5_18} and \eqref{eq5_19}. This implies that there exist arbitrarily large $\tau$ satisfying \eqref{eqT} and \eqref{eqQ}.
    
\end{proof}

Thus by putting $x=\frac{\log T}{3(\log\log T)^3}$  in \eqref{eq5_17} and using the fact that there is a critical point $\rho'\in\mathscr{C}_1+\di\tau$ and $T\leq \tau \leq 2T$, we get

$$\limsup\limits_{\gamma' \rightarrow\infty}\frac{|L(\rho',\chi)|}{\log \log\gamma'}\geq e^{C_0}\epsilon\frac{\phi(q)}{q}.$$
By taking $\epsilon$ arbitrary close to 1, we can get
$$\limsup\limits_{\gamma' \rightarrow\infty}\frac{|L(\rho',\chi)|}{\log \log\gamma'}\geq e^{C_0}\frac{\phi(q)}{q}.$$
\end{proof}
\section{Proof of Theorem 3}
\begin{proof}
From Lemma \ref{lemma1_1},
\begin{equation}\label{est3}
  \sum_{1< n \leq x}\frac{\chi(n)\Lambda(n)}{n^s\log n}=\sum_{p\leq x}\log\bigg(1-\frac{\chi(p)}{p^s}\bigg)^{-1}+O\bigg(\frac{1}{\sqrt x\log x}\bigg).  
\end{equation}

Since $\Re \log(1-\chi(p)p^{-s})^{-1}\geq -\log(1+1/p)$, when $|\chi(p)|=1$, we get
\begin{align*}
    \Re\sum_{1< n \leq x}\frac{\chi(n)\Lambda(n)}{n^s\log n} &\geq -\sum_{p\leq x}\log\bigg(1+\frac{1}{p}\bigg)+\sum_{p\mid q}\log \bigg(\frac{p+1}{p}\bigg) +O\bigg(\frac{1}{\sqrt x\log x}\bigg).\\
    \end{align*}
Now since,
$$-\sum_{p\leq x}\log\bigg(1+\frac{1}{p}\bigg)=-\sum_{p\leq x}\log\bigg(1-\frac{1}{p^2}\bigg)+\sum_{p\leq x}\log\bigg(1-\frac{1}{p}\bigg),$$
and
$$-\sum_{p\leq x}\log\bigg(1-\frac{1}{p^2}\bigg)=\log(\pi^2/6)+O(1/x),$$
and
$$\sum_{p\leq x}\log\bigg(1-\frac{1}{p}\bigg)=-\log \log x-C_0+O\bigg(\frac{1}{\log x}\bigg),$$
we have, 
\begin{equation}\label{eq21_3}
    \Re \sum_{1< n \leq x}\frac{\chi(n)\Lambda(n)}{n^s\log n}\geq  -\log \log x-C_0+\log(\pi^2/6)+\sum_{p\mid q}\log\bigg(\frac{p+1}{p}\bigg)+ O\bigg(\frac{1}{\log x}\bigg).
\end{equation}

By Lemma \ref{lemma4_1}, we have
\begin{equation}\label{eq22_3}
  \log|L(s,\chi)|=\Re \sum_{1< n \leq \log^2 T}\frac{\chi(n)\Lambda(n)}{n^s\log n}+O\bigg(\frac{1}{\log \log T}\bigg).  
\end{equation}
Combining \eqref{eq21_3} and \eqref{eq22_3} and putting $x=\log^2 T$ with $4\leq T \leq \gamma'\leq 2T$, we get

$$\liminf\limits_{\gamma' \rightarrow\infty}|L(\rho',\chi)|\log\log \gamma'\geq \frac{\pi^2e^{-C_0}}{12}\prod_{p\mid q}\frac{p+1}{p}.$$
\end{proof}
Remark: In the proof of Theorem 3 we did not use any special property of critical points of the zeta function. One can also recover
the following classical bound.\\
Let $s=\sigma +\di t$ be a point with $\sigma\geq1$, then
$$\liminf\limits_{t\rightarrow\infty}|L(s,\chi)|\log\log t\geq \frac{\pi^2e^{-C_0}}{12}\prod_{p\mid q}\frac{p+1}{p}.$$
\section{\textbf{Lemmas related to the proof of Theorem 4}}
Let $\epsilon,\delta,q$ be as defined before. For $x>q^{1/\epsilon}$, define a function $$Z_x(s)=\sum_{n\leq x}\frac{\Lambda(n)b'(n)}{n^s},$$
where $b'(n)$ is a totally multiplicative function such that \\
$$
b'(p)=
\begin{cases}
-1 &\text {if p$\mid$ q},\\
-\bar\chi(p) &\text { $p\leq x^\epsilon$ and p$\nmid$ q},\\
1 & \text{if $x^\epsilon<p\leq mx^\epsilon$},\\
-1 & \text{if $mx^\epsilon<p\leq x^\delta$},\\
1& \text{if $x^\delta<p$}.
\end{cases}
$$

Define 
\begin{equation}\label{eqs2}
  S_2=\frac{L'}{L}(1,\bar\chi)+2\sum_{p=2}^{\infty}\frac{\bar\chi(p)^2\log p}{p^2-\bar\chi(p)^2}-\sum_{p\mid q}\frac{\log p}{p+1}.  
\end{equation}

Let $m$ be such that $$2\log m=-\Re S_2+|S_2|+1.$$

 \begin{lemma}\label{lemma1_4}
If $x^\epsilon>q$, $\Re s\geq 1$ and $$|s-1|\leq\frac{6|S_2|+3}{|2\delta^2-1-\epsilon^2|\log^2 x},$$ then $$Z_x(s)=\frac{(s-1)}{2}(2\delta^2-1-\epsilon^2)\log^2 x+S_2+2\log m+O\bigg(\frac{1}{\log x}\bigg).$$
 Furthermore, $Z_x(s) $ has a root at $$s=1+\frac{2(S_2+2\log m)}{(-2\delta^2+1+\epsilon^2)\log^2 x}+O\bigg(\frac{1}{\log^3 x}\bigg).$$
\end{lemma}
\begin{proof}
We will split $Z_x(s)$ as follows. 
$$Z_x(s)=\sum_{n\leq x^\epsilon}\frac{\Lambda(n)b'(n)}{n^s}+\sum_{ x^\epsilon < n \leq mx^\epsilon}\frac{\Lambda(n)b'(n)}{n^s}+\sum_{ mx^\epsilon < n \leq x^\delta}\frac{\Lambda(n)b'(n)}{n^s}+\sum_{x^\delta< n\leq x}\frac{\Lambda(n)b'(n)}{n^s}.$$
We represent the four sums by $T_1, T_2, T_3,T_4$ respectively.

We estimate these sums like we did in part 2. First we estimate $T_1.$
Using \eqref{eqap} and \eqref{eqs2}, we get
\begin{equation}\label{eq84}
    \begin{split}
    T_1 &=\sum_{n\leq x^\epsilon}\frac{\Lambda(n)b'(n)}{n^s}=\frac{L'}{L}(s,\bar\chi)+2\sum_{p^{2k}\leq x^\epsilon}\frac{\bar\chi(p)^{2k}\log p}{p^{2ks}}\\
    &+ \sum_{p^k\leq x^\epsilon\atop p|q}\frac{(-1)^k\log p}{p^{ks}}+O\big(x^{\epsilon(1-\sigma)}\exp(-\tilde c\sqrt{\log x})\big).\\
 \end{split}
\end{equation}

We have the following approximations.
\begin{equation}\label{eq85}
\begin{split}
\sum_{p^{2k}\leq x^\epsilon}\frac{\bar\chi(p)^{2k}\log p}{p^{2ks}}&= \sum_{p^{2k}\leq x^\epsilon}\frac{\bar\chi(p)^{2k}\log p}{p^{2k}}+O(1/\log x)\\
&=\sum_{k=1}^{\infty}\sum_{p=2}^{\infty}\frac{\bar\chi(p)^{2k}\log p}{p^{2k}}+O(1/\log x)\\
&=\sum_{p=2}^{\infty}\frac{\bar\chi(p)^2\log p}{p^2-\bar\chi(p)^2}+O(1/\log x)
\end{split}
\end{equation}
and
\begin{equation}\label{eq86}
    \begin{split}
        \sum_{p^k\leq x^\epsilon\atop p|q}\frac{(-1)^k\log p}{p^{ks}}&=\sum_{p^k\leq x^\epsilon\atop p|q}\frac{(-1)^k\log p}{p^{k}}+O(1/\log x)\\
        &=\sum_{k=1}^{\infty}\sum_{ p|q}\frac{(-1)^k\log p}{p^{k}}+O(1/\log x)\\
        &=-\sum_{p\mid q}\frac{\log p}{p+1}+O(1/\log x).
\end{split}
\end{equation}

and from \eqref{newbound}
$$
\frac{L'}{L}(s,\bar\chi)-\frac{L'}{L}(1,\bar \chi)\ll \frac{1}{\log x}.
$$
Combining \eqref{eq84},\eqref{eq85} and \eqref{eq86}, we get
\be
\begin{split}
    T_1 &=\frac{L'}{L}(1,\bar\chi)+2\sum_{p^{2k}\leq x^\epsilon}\frac{\bar\chi(p)^{2k}\log p}{p^{2k}}+\sum_{p^k\leq x^\epsilon \atop p|q}\frac{(-1)^k\log p}{p^{k}}+O\big(1/\log x\big)\\
    &=\frac{L'}{L}(1,\bar\chi)+2\sum_{p=2}^{\infty}\frac{\bar\chi(p)^2\log p}{p^2-\bar\chi(p)^2}-\sum_{p\mid q}\frac{\log p}{p+1}+O(1/\log x).\\
\end{split}
\ee
Recall that, 
\begin{equation}
  S_2=\frac{L'}{L}(1,\bar\chi)+2\sum_{p=2}^{\infty}\frac{\bar\chi(p)^2\log p}{p^2-\bar\chi(p)^2}-\sum_{p\mid q}\frac{\log p}{p+1}.  
\end{equation}
So we have
\begin{equation}\label{eq87n}
    T_1=S_2+O\big(1/\log x\big).
\end{equation}

We estimate $T_2$ as follows.
\begin{equation}\label{eq88n}
\begin{split}
T_2=\sum_{x^\epsilon< n\leq mx^\epsilon}\frac{\Lambda(n)b'(n)}{n^s}&=\sum_{x^\epsilon< p\leq mx^\epsilon}\frac{\log p}{p^s}+O\bigg(\sum_{k=2 }^{\log (mx)^\epsilon/\log 2}\sum_{p=x^{\epsilon/k}}^{(mx)^{\epsilon/k}}\frac{\log p}{p^k}\bigg)\\ &=\sum_{x^\epsilon< p\leq mx^\epsilon}\frac{\log p}{p^s}+O(\log^{-1} x)\\
&=\sum_{x^\epsilon< p\leq mx^\epsilon}\frac{\log p}{p}+O(\log^{-1} x)\\ 
&= \log m +O(\log^{-1} x).  
\end{split}    
\end{equation}

$T_3$ is estimated as follows.
\begin{equation}
\bigg|T_3+\sum_{mx^\epsilon <n\leq x^\delta}\frac{\Lambda(n)}{n^s}\bigg|\ll \sum_{k=2}^{\delta \log x/\log 2}\sum_{p=(mx^{\epsilon})^{1/k}}^{x^{\delta/k}}\frac{\log p}{p^k}\ll \sum_{k=2}^{\log x/\log 2}\frac{1}{(mx)^{\epsilon(k-1)/k}}\ll(\log x)^{-1}.    
\end{equation}
So we get 
\be\label{eq89n}
T_3=V_{mx^\epsilon}(s)-V_{x^\delta}(s).
\ee

Similarly we can prove
\begin{equation}
\bigg|T_4-\sum_{x^\delta <n\leq x}\frac{\Lambda(n)}{n^s}\bigg|\ll \sum_{k=2 }^{\log x/\log 2}\sum_{p=x^{\delta/k}}^{x^{1/k}}\frac{\log p}{p^k}\ll \sum_{k=2 }^{\log x/\log 2}\frac{1}{x^{\delta (k-1)/k}}\ll(\log x)^{-1}.
\end{equation}
So we have 
\be \label{eq811n}
T_4=V_{x}(s)-V_{x^\delta}(s).
\ee

Combining \eqref{eq87n}, \eqref{eq88n}, \eqref{eq89n} and \eqref{eq811n}, we obtain
\begin{equation}\label{zequation}
\begin{split}
Z_x(s)&=-2V_{x^\delta}(s)+V_x(s)+V_{mx^\epsilon}(s)+S_2+\log m+O(1/\log x)\\
&=-\frac{2x^{\delta(1-s)}}{1-s}+\frac{x^{1-s}}{1-s}+\frac{(mx^{\epsilon})^{(1-s)}}{1-s}+S_2+\log m+O(1/\log x)\\
&=\bigg(\frac{-2e^{\delta(1-s)\log x}+e^{(1-s)\log x}+e^{(\log m+\epsilon\log x)(1-s)}}{1-s}\bigg)+S_2+\log m+O(1/\log x)\\
&=\frac{(s-1)}{2}\big[(2\delta^2-1-\epsilon^2)\log^2 x-(\log^2m+2\epsilon\log m \log x)\big]+S_2+2\log m+O({1/\log x}).
\end{split}
\end{equation}

Now since $$|s-1|\leq\frac{6|S_2|+3}{|2\delta^2-1-\epsilon^2|\log^2 x}$$
we have,
$$Z_x(s)= \frac{(s-1)}{2}(2\delta^2-1-\epsilon^2)\log^2 x+S_2+2\log m+O\bigg(\frac{1}{\log x}\bigg)$$

and 
$Z_x(s) $ has a root at $$s=1+\frac{2(S_2+2\log m)}{(-2\delta^2+1+\epsilon^2)\log^2 x}+O\bigg(\frac{1}{\log^3 x}\bigg).$$
For the given choice of $m$, we can see that $S_2+2\log m>0$.
Since $2(S_2+2\log m)/(-2\delta^2+1+\epsilon^2)>0$, we can say that the root has real part greater than 1.
\end{proof}
\begin{center}
We now proceed to prove Theorem 4.        
\end{center}

\section{\textbf{The proof of Theorem 4}}
\begin{proof}
Define circles $\mathscr{C}_0$ and $\mathscr{C}_1$ as
$$\mathscr{C}_0= \bigg\{s=1+\frac{2(S_2+2\log m)}{(-2\delta^2+1+\epsilon^2)\log^2 x}+\frac{(\Re S_2+2\log m)}{(-2\delta^2+1+\epsilon^2)\log^2 x}e^{\di\phi}:0\leq \phi\leq 2\pi\bigg\},$$
$$\mathscr{C}_1= \bigg\{s=1+\frac{2(S_2+2\log m)}{(-2\delta^2+1+\epsilon^2)\log^2 x}+\frac{(\Re S_2+2\log m)}{2(-2\delta^2+1+\epsilon^2)\log^2 x}e^{\di\phi}:0\leq \phi\leq 2\pi\bigg\}.$$
Let $\tau$ be a number such that 
\begin{equation}
\|\tau\log p/2\pi-\arg{\bar b'(p)/2\pi}\|\leq 1/\log^2 x,   
\end{equation}
for all $p\leq x$.\\
For such $\tau$ and for all $p\leq x$, we get 

\begin{equation}\label{eq1_4}
   |p^{-\di\tau}-b'(p)|\ll 1/\log^2 x. 
\end{equation}

For a $\tau$ satisfying \eqref{eq1_4} and for $s\in \mathscr{C}_1$, we have
\begin{equation}
    \begin{split}
    |V_x(s+\di\tau)-Z_x(s)|\ll\sum_{ n \leq x}\bigg|\frac{\Lambda(n)}{n^{s+\di\tau}}-\frac{\Lambda(n)b'(n)}{n^s}\bigg|\ll\sum_{p^k \leq x}\frac{k\log p}{p^k}|p^{-\di\tau}-b'(p)|\ll\frac{1}{\log x}.
    \end{split}
\end{equation}

So we have
\begin{equation}
    V_x(s+\di\tau)=Z_x(s)+O\bigg(\frac{1}{\log x}\bigg).
\end{equation}
For a sufficiently large $x$, define $M_x(s)$ as $$M_x(s)=\sum_{1<n\leq x}\frac{\Lambda(n)c'(n)}{n^s\log n},$$
where $c'(n)$ is a completely multiplicative function such that 
$$
c'(p)=
\begin{cases}
- 1 &\text{if\qquad} p\mid q,\\
 \chi(p)b'(p) &\text{if\qquad}p \nmid q.\\
 \end{cases}
$$

For $s\in \mathscr{C}_1$, by Lemma \ref{lemma1_1}, we have,
\begin{equation}
\begin{split}
\sum_{1<n\leq x \atop (q,n)>1}\frac{\Lambda(n)c'(n)}{n^s\log n}&=\sum_{p\mid q}\log\bigg(1-\frac{c'(p)}{p^s}\bigg)^{-1}+O\bigg(\frac{1}{\sqrt{x}\log x}\bigg)    \\
&=-\sum_{p\mid q}\log\bigg(1+\frac{1}{p^s}\bigg)+O\bigg(\frac{1}{\sqrt{x}\log x}\bigg)  \\
&=-\sum_{p\mid q}\log\bigg(1+\frac{1}{p}\bigg)+O\bigg(\frac{1}{\log^2 x}\bigg).
\end{split}
\end{equation}
We get the last equality because
$$
\log\bigg(1+\frac{1}{p^s}\bigg)-\log\bigg(1+\frac{1}{p}\bigg)\ll\frac{1}{\log^2 x}.
$$ 

Thus, for $s\in\mathscr{C}_1$ we get
\begin{equation}
  \begin{split}
     \sum_{1<n\leq x}\frac{\Lambda(n)\chi(n)}{n^{s+\di\tau}\log n}-M_x(s)- \sum_{p\mid q}\log \frac{p+1}{p} &=\sum_{1<n\leq x}\frac{\Lambda(n)\chi(n)}{n^{s+\di\tau}\log n}-\sum_{1<n\leq x}\frac{\Lambda(n)c'(n)}{n^s\log n}+\sum_{p\mid q}\sum_{m=1}^{\infty}\frac{(-1)^{m}}{mp^m} \\
     &=\sum_{1<n\leq x}\frac{\Lambda(n)\chi(n)}{n^{s+\di\tau}\log n}-\sum_{1<n\leq x \atop (q,n)=1}\frac{\Lambda(n)c'(n)}{n^s\log n}+O\bigg(\frac{1}{\log x}\bigg)\\
     &\leq\sum_{1<n\leq x \atop (q,n)=1}\frac{\Lambda(n)}{n\log n}\bigg|\frac{\chi(n)}{n^{\di\tau}}-c'(n)\bigg|+O\bigg(\frac{1}{\log x}\bigg)\\
     &\ll\sum_{1<p^k\leq x \atop p \nmid q}\frac{1}{p^k} \bigg|\frac{\chi(p)}{p^{\di\tau}}-c'(p)\bigg|+O\bigg(\frac{1}{\log x}\bigg).\\
\end{split}
  \end{equation}
We know that for a $\tau$ satisfying \eqref{eq1_4} and $p\leq x$, we have
$$|p^{-\di \tau}-b'(p)|\ll 1/\log^2 x.$$
So for $p\leq x$ and $p\nmid q$, we have
\begin{equation*}
    \begin{split}
     |p^{-\di \tau}-b'(p)|=|p^{-\di \tau}\chi(p)-b'(p)\chi(p)|=\bigg|\frac{\chi(p)}{p^{\di\tau}}-c'(p)\bigg|\ll 1/\log^2 x.   
    \end{split}
\end{equation*}
Thus,
$$\sum_{1<p^k\leq x \atop p \nmid q}\frac{1}{p^k} \bigg|\frac{\chi(p)}{p^{\di\tau}}-c'(p)\bigg|\ll \frac{1}{\log x}.$$
Combining all the above estimates, we have
\begin{equation}\label{eq22_4}
    \sum_{1<n\leq x}\frac{\Lambda(n)\chi(n)}{n^{s+\di\tau}\log n}=M_x(s)+ \sum_{p\mid q}\log \frac{p+1}{p}+O\bigg(\frac{1}{\log x}\bigg).
\end{equation}

Since
$$M'_x(s)\ll \sum_{1\leq n \leq x}\frac{\Lambda(n)c'(n)}{n}\ll \log x,$$
by Liouville's theorem on $\mathscr{C}_1$, we get 
\begin{equation}\label{eq23_4}
  M_x(s)=M_x(1)+O\bigg(\frac{1}{\log x}\bigg).  
\end{equation}
We have  $$M_x(1)=\sum_{1<n\leq x^\epsilon}\frac{\Lambda(n)c'(n)}{n\log n}+\sum_{x^\epsilon<n\leq mx^\epsilon}\frac{\Lambda(n)c'(n)}{n\log n}+\sum_{mx^\epsilon<n\leq x^\delta}\frac{\Lambda(n)c'(n)}{n\log n}+\sum_{x^\delta<n\leq x}\frac{\Lambda(n)c'(n)}{n\log n}.$$

Using the similar method from Part 2, we can see that the second, third and fourth terms on the right-hand side of the above equation are $O(1/\log x).$ Thus, we get the following expression for $M_x(1)$.

\begin{equation}\label{eq24_4}
    \begin{split}
M_x(1)=\sum_{1<n\leq x}\frac{\Lambda(n)c'(n)}{n\log n}&=\sum_{1<n\leq x^\epsilon}\frac{\Lambda(n)c'(n)}{n\log n}+O\bigg(\frac{1}{\log x}\bigg)\\&=\sum_{p\leq x^\epsilon}\log\bigg(1-\frac{c'(p)}{p}\bigg)^{-1}+O\bigg(\frac{1}{\log x}\bigg)\\
&= -\sum_{p\leq x^\epsilon}\log\bigg(1+\frac{1}{p}\bigg) +O\bigg(\frac{1}{\log x}\bigg)\\
    &=-\log \log x^\epsilon-C_0+\log(\pi^2/6)+ O\bigg(\frac{1}{\log x}\bigg) .
\end{split}
\end{equation}
We get the last line in the above equation from Theorem 2.7 of \cite{classical}.\\
Combining \eqref{eq22_4},\eqref{eq23_4} and \eqref{eq24_4}, we get

\begin{equation}\label{eq25_4}
 \Re \sum_{1<n\leq x}\frac{\Lambda(n)\chi(n)}{n^{s+\di\tau}\log n}\leq -\log\log x -C_0+\log(\pi^2/6)+\sum_{p\mid q}\log \frac{p+1}{p} -\log \epsilon+O\bigg(\frac{1}{\log x}\bigg).   
\end{equation}

With $T_x(s)$, $Q_x(s)$ and $Q'_x(s)$ as defined in part 2. For $\Re s>1$, we have
$$\log L(s,\chi)= \sum_{1<n\leq x}\frac{\Lambda(n)\chi(n)}{n^{s}\log n}+T_x(s)+O\bigg(\frac{1}{\log x}\bigg),$$ and
$$-\frac{\zeta'}{\zeta}(s)=\sum_{1\leq n\leq x}\frac{\Lambda(n)}{n^s}-Q'_x(s)+O\bigg(\frac{1}{\log x}\bigg).$$
If $\tau$ satisfies \eqref{eqT} and \eqref{eqQ}, then we get
\begin{equation}\label{eq914n}
    \Re \log L(s+\di\tau,\chi)\leq -\log\log x -C_0+\log(\pi^2/6)+\sum_{p\mid q}\log \frac{p+1}{p} -\log \epsilon+O\bigg(\frac{1}{\log x}\bigg),
\end{equation}
and
\begin{equation}\label{eq29_4}
  \begin{split}
  -\frac{\zeta'}{\zeta}(s+\di\tau)=Z_x(s)+O\bigg(\frac{1}{\log x}\bigg).   \end{split}
  \end{equation}
From \eqref{eq29_4}, we have
\begin{equation}
-\frac{\zeta'}{\zeta}(s+\di\tau)-Z_x(s)\ll\frac{1}{\log x}.   
\end{equation}
If $Z_x(s)\asymp 1$ on $\mathscr{C}_1$ then
\begin{equation}\label{rouche12}
-\frac{\zeta'}{\zeta}(s+\di\tau)-Z_x(s)\ll\frac{1}{\log x}\ll Z_x(s).   
\end{equation}
If $Z_x(s)$ has a root for $s\in\mathscr{C}_1$ and $Z_x(s)\asymp 1$ on the boundary of $\mathscr{C}_1$, then by Rouche's theorem $\zeta'(s)$ will have a zero $\rho'$ for $\rho'\in\mathscr{C}_1+\di \tau$ and from \eqref{eq914n}, we get
\begin{equation}\label{eq918nn}
 |L(\rho',\chi)|\leq \frac{1}{\epsilon\log x} e^{-C_0} \frac{\pi^2}{6} \prod_{p|q}\frac{p+1}{p}e^{O(1/\log x)}.
\end{equation}
 
Finally, to prove Theorem 4, we want the following to be true:
\begin{itemize}
  \item $Z_x(s)$ to have a root inside $\mathscr{C}_1$
  \item $Z_x(s)\asymp 1$ on $\mathscr{C}_1.$ 
  \item Existence of arbitrarily large $\tau$ satisfying \eqref{eqT} and \eqref{eqQ}.
\end{itemize}
The fulfillment of the first two conditions above depends on making the correct choice for $m$. Recall that 
\begin{equation}\label{mvalue1}
  2\log m=-\Re S_2+|S_2|+1.  
\end{equation}

\begin{lemma}
    $Z_x(s)$ has a root inside $\mathscr{C}_1$.
\end{lemma}
\begin{proof}
    To prove this, we need to check that the distance of the root of $Z_x(s)$ from the center of $\mathscr{C}_1$ is less than the radius.

From Lemma \ref{lemma1_4}, we know that $Z_x(s)$ has a root at 
$$s=1+\frac{2(S_2+2\log m)}{(-2\delta^2+1+\epsilon^2)\log^2 x}+O\bigg(\frac{1}{\log^3 x}\bigg).$$
Its distance from the center $\mathscr{C}_1$ is $O(1/\log^3 x)$. The radius of $\mathscr{C}_1$ is $$\frac{(\Re S_2+2\log m)}{(-2\delta^2+1+\epsilon^2)\log^2 x}.$$ Thus proving the assertion.
\end{proof}

\begin{lemma}
   $Z_x(s)\asymp 1$ on $\mathscr{C}_1.$ 
\end{lemma}
\begin{proof}
    To check how $Z_x(s)$ varies on $\mathscr{C}_1$, we plug the equation of $\mathscr{C}_1$ in the equation of $Z_x(s)$. We get,
\begin{equation*}
\begin{split}
Z_x(\mathscr{C}_1)&=\frac{1}{2}\bigg(\frac{2(S_2+2\log m)}{(-2\delta^2+1+\epsilon^2)\log^2 x}+\frac{(\Re S_2+2\log m)}{2(-2\delta^2+1+\epsilon^2)\log^2 x}e^{\di\phi}\bigg)(2\delta^2-1-\epsilon^2)\log^2x+S_2+2\log m\\
&+O\bigg(\frac{1}{\log x}\bigg)\\
&=\frac{1}{4}(\Re S_2+2\log m)e^{\di\phi}+O\bigg(\frac{1}{\log x}\bigg)\asymp 1.
\end{split}
\end{equation*}
    Thus, $Z_x(s)\asymp 1$ on the boundary of $\mathscr{C}_1.$

\end{proof}
 
The existence of arbitrarily large $\tau$ satisfying \eqref{eqT} and \eqref{eqQ}is proved in the proof of Theorem 2 .
Using $x={\log T}/{3(\log \log T)^3}$ in \eqref{eq918nn} and the existence of a critical point $\rho'\in \mathscr{C}_1+\di\tau$ for $T\leq \tau \leq 2T$, we get
$$\liminf\limits_{\gamma' \rightarrow\infty}|L(\rho',\chi)|\log\log \gamma'\leq \frac{\pi^2e^{-C_0}}{6\epsilon}\prod_{p\mid q}\frac{p+1}{p}.$$
By taking $\epsilon$ arbitrarily close to 1,$$\liminf\limits_{\gamma' \rightarrow\infty}|L(\rho',\chi)|\log\log \gamma'\leq \frac{\pi^2e^{-C_0}}{6}\prod_{p\mid q}\frac{p+1}{p}.$$
\end{proof}

\section{Remarks}
One can check that even if we replace the the critical points of the Riemann zeta function with $1+\di t$, the bounds in Theorems 2 and 4 still remain the same. By Kronecker's theorem we see that when $x=\log T/3(\log\log T)^3$, we can choose $t$ such that ${\chi(n)}{n^{-\di t}}\approx 1$ for all $n\leq x$, thus getting the same extreme values.

\section{Acknowledgement}
The author gives sincere thanks to his doctoral advisor Steven M. Gonek for
introducing the problem in this paper and also for providing guidance and support
during the process of its study. Professor Gonek also read an earlier version of
this paper and made many useful suggestions which significantly improved the exposition. The author would also like to thanks professor Atul Dixit and professor Bruce Berndt for their valuable suggetions.


\begin{thebibliography}{9}

\bibitem{papernew1}
C. Aistleitner, K. Mahatab, M. Munsch, Extreme Values of the Riemann Zeta Function on the 1-Line, \textit{International Mathematical Research Notices} (2019) no.22, 6924–6932.

\bibitem{papernew2} C. Aistleitner, K. Mahatab, M. Munsch, A. Peyrot, On large values of $L(\sigma,\chi)$, \textit{Q.J.Math.}\textbf{70} (2019), no.3, 831-848.

\bibitem{paper5}
E. C. Tichmarsh, On the function $1/\zeta(1+\di t)$,\textit{Quart. J. math.(Oxford)},\textbf{4} (1933),64-70.

\bibitem{Titch}  E. C. Titchmarsh, The theory of the Riemann
zeta-function (2nd edition),  Oxford University Press, Oxford, 1986.

\bibitem{paper2}
H. L. Montgomery and R. C. Vaughan,Hilbert's inequality, \textit{J. London Math Soc.}(2) \textbf{8} (1974), 73-82.

\bibitem{classical} 
	H. L. Montgomery and R. C. Vaughan,
	\textit{Multiplicative Number Theory 1: Classical Theory}, 
	Cambridge University Press, Cambridge,2007.

\bibitem{paper4}
J. E. Littlewood, On the function $1/\zeta(1+\di t)$, \textit{Proc. London math. Soc.}(2) \textbf{27} (1928), 349-357.

\bibitem{paper3}
J. E. Littlewood, On the Riemann zeta function, \textit{Proc. London math. Soc.}(2) \textbf{24} (1926), 175-201.

\bibitem{Mont-Thomp}
H. L. Montgomery and  J. G. Thompson,   Geometric properties of the zeta function \textit{Acta Arith.} 
\textbf{155} (2012), 373–96.

\bibitem{papernew5}
S. Chorge, Extreme values of the Riemann zeta function at its critical points,
\textit{https://arxiv.org/abs/2110.14229}(accepted by Acta Arithmetica but yet to published)


\bibitem{paper}
S. M. Gonek and H. L. Montgomery,Extreme values of zeta function at critical points, \textit{Q. J. Math.} \textbf{67} (2016), no. 3, 483–505.

\bibitem{papernew3}
J. Li, Large values of Dirichlet $L$-functions at zeros of a class of $L$-functions, \textit{Canad. J. Math.}\textbf{73} (2021), no. 6, 1459-1505.

\bibitem{papernew4}
Y. Lamzouri, X. Li, K. Soundararajan, Conditional bounds for the least quadratic non-residue and related problems, \textit{Math.Comp.}\textbf{84} (2015), no. 295, 2391-2412.
\end{thebibliography}
\end{document}